\documentclass[11pt,a4paper]{article}
\usepackage{amsmath, amssymb}
\usepackage{latexsym, color}
\usepackage{epsfig}

\numberwithin{equation}{section}

\begin{document}

\newtheorem{thm}{Theorem}[section]
\newtheorem{cor}[thm]{Corollary}
\newtheorem{lem}[thm]{Lemma}
\newtheorem{prop}[thm]{Proposition}
\newtheorem{definition}[thm]{Definition}
\newtheorem{rem}[thm]{Remark}
\newtheorem{Ex}[thm]{EXAMPLE}
\def\nm{\noalign{\medskip}}

\bibliographystyle{plain}

%\numberwithin{equation}{section}

\newcommand{\qed}{\hfill \ensuremath{\square}}
\newcommand{\ds}{\displaystyle}
\newcommand{\pf}{\noindent {\sl Proof}. \ }
\newcommand{\p}{\partial}
\newcommand{\pd}[2]{\frac {\p #1}{\p #2}}
\newcommand{\norm}[1]{\| #1 \|}
\newcommand{\dbar}{\overline \p}
\newcommand{\eqnref}[1]{(\ref {#1})}
\newcommand{\na}{\nabla}
\newcommand{\one}[1]{#1^{(1)}}
\newcommand{\two}[1]{#1^{(2)}}

\newcommand{\Abb}{\mathbb{A}}
\newcommand{\Cbb}{\mathbb{C}}
\newcommand{\Ibb}{\mathbb{I}}
\newcommand{\Nbb}{\mathbb{N}}
\newcommand{\Kbb}{\mathbb{K}}
\newcommand{\Rbb}{\mathbb{R}}
\newcommand{\Sbb}{\mathbb{S}}

\renewcommand{\div}{\mbox{div}~}

\newcommand{\la}{\langle}
\newcommand{\ra}{\rangle}

\newcommand{\Hcal}{\mathcal{H}}
\newcommand{\Ical}{\mathcal{I}}
\newcommand{\Lcal}{\mathcal{L}}
\newcommand{\Kcal}{\mathcal{K}}
\newcommand{\Dcal}{\mathcal{D}}
\newcommand{\Pcal}{\mathcal{P}}
\newcommand{\Qcal}{\mathcal{Q}}
\newcommand{\Scal}{\mathcal{S}}

%%%%%%%%%%%%%%%%%%%%%%%%%%%%%%%%%%%%%%%%%%%%
% define bold face
%%%%%%%%%%%%%%%%%%%%%%%%%%%%%%%%%%%%%%%%%%%%
\def\Ba{{\bf a}}
\def\Bb{{\bf b}}
\def\Bc{{\bf c}}
\def\Bd{{\bf d}}
\def\Be{{\bf e}}
\def\Bf{{\bf f}}
\def\Bg{{\bf g}}
\def\Bh{{\bf h}}
\def\Bi{{\bf i}}
\def\Bj{{\bf j}}
\def\Bk{{\bf k}}
\def\Bl{{\bf l}}
\def\Bm{{\bf m}}
\def\Bn{{\bf n}}
\def\Bo{{\bf o}}
\def\Bp{{\bf p}}
\def\Bq{{\bf q}}
\def\Br{{\bf r}}
\def\Bs{{\bf s}}
\def\Bt{{\bf t}}
\def\Bu{{\bf u}}
\def\Bv{{\bf v}}
\def\Bw{{\bf w}}
\def\Bx{{\bf x}}
\def\By{{\bf y}}
\def\Bz{{\bf z}}
\def\BA{{\bf A}}
\def\BB{{\bf B}}
\def\BC{{\bf C}}
\def\BD{{\bf D}}
\def\BE{{\bf E}}
\def\BF{{\bf F}}
\def\BG{{\bf G}}
\def\BH{{\bf H}}
\def\BI{{\bf I}}
\def\BJ{{\bf J}}
\def\BK{{\bf K}}
\def\BL{{\bf L}}
\def\BM{{\bf M}}
\def\BN{{\bf N}}
\def\BO{{\bf O}}
\def\BP{{\bf P}}
\def\BQ{{\bf Q}}
\def\BR{{\bf R}}
\def\BS{{\bf S}}
\def\BT{{\bf T}}
\def\BU{{\bf U}}
\def\BV{{\bf V}}
\def\BW{{\bf W}}
\def\BX{{\bf X}}
\def\BY{{\bf Y}}
\def\BZ{{\bf Z}}

%%%%%%%%%%%%%%%%%%%%%%%%%%%%%%%%%%%%%%%%%%%%%%%%%%%%
% Abbreviate definitions of greek symbols
%%%%%%%%%%%%%%%%%%%%%%%%%%%%%%%%%%%%%%%%%%%%%%%%%%%%

\newcommand{\Ga}{\alpha}
\newcommand{\Gb}{\beta}
\newcommand{\Gd}{\delta}
\newcommand{\Ge}{\epsilon}
\newcommand{\Gve}{\varepsilon}
\newcommand{\Gf}{\phi}
\newcommand{\Gvf}{\varphi}
\newcommand{\Gg}{\gamma}
\newcommand{\Gc}{\chi}
\newcommand{\Gi}{\iota}
\newcommand{\Gk}{\kappa}
\newcommand{\Gvk}{\varkappa}
\newcommand{\Gl}{\lambda}
\newcommand{\Gn}{\eta}
\newcommand{\Gm}{\mu}
\newcommand{\Gv}{\nu}
\newcommand{\Gp}{\pi}
\newcommand{\Gt}{\theta}
\newcommand{\Gvt}{\vartheta}
\newcommand{\Gr}{\rho}
\newcommand{\Gvr}{\varrho}
\newcommand{\Gs}{\sigma}
\newcommand{\Gvs}{\varsigma}
\newcommand{\Gj}{\tau}
\newcommand{\Gu}{\upsilon}
\newcommand{\Go}{\omega}
\newcommand{\Gx}{\xi}
\newcommand{\Gy}{\psi}
\newcommand{\Gz}{\zeta}
\newcommand{\GD}{\Delta}
\newcommand{\GF}{\Phi}
\newcommand{\GG}{\Gamma}
\newcommand{\GL}{\Lambda}
\newcommand{\GP}{\Pi}
\newcommand{\GT}{\Theta}
\newcommand{\GS}{\Sigma}
\newcommand{\GU}{\Upsilon}
\newcommand{\GO}{\Omega}
\newcommand{\GX}{\Xi}
\newcommand{\GY}{\Psi}

%%%%%%%%%%
\newcommand{\beq}{\begin{equation}}
\newcommand{\eeq}{\end{equation}}

\title{Optimal estimates and asymptotics for the stress concentration between closely located stiff inclusions\thanks{\footnotesize This work is supported by Korean Ministry of Education, Sciences and Technology through NRF grants Nos. 2010-0017532, 2013R1A1A1A05009699 and  2011-0009671.}  }

\author{Hyeonbae Kang\thanks{Department of Mathematics, Inha University, Incheon
402-751, Korea (hbkang@inha.ac.kr, hdlee@inha.ac.kr).}  \and Hyundae Lee\footnotemark[2] \and KiHyun Yun\thanks{Department of Mathematics, Hankuk University of Foreign Studies,
Yongin-si, Gyeonggi-do 449-791, Korea (gundam@hufs.ac.kr).}}

\date{\today}

\maketitle

\begin{abstract}
If stiff inclusions are closely located, then the stress, which is the gradient of the solution, may become arbitrarily large as the distance between two inclusions tends to zero. In this paper we investigate the asymptotic behavior of the stress concentration factor, which is the normalized magnitude of the stress concentration, as the distance between two inclusions tends to zero. For that purpose we show that the gradient of the solution to the case when two inclusions are touching decays exponentially fast near the touching point.  We also prove a similar result when two inclusions are closely located and there is no potential difference on boundaries of two inclusions. We then use these facts to show that the stress concentration factor converges to a certain integral of the solution to the touching case as the distance between two inclusions tends to zero. We then present an efficient way to compute this integral.
\end{abstract}

\noindent {\footnotesize {\bf AMS subject classifications.} 35J25, 73C40}

\noindent {\footnotesize {\bf Key words.} conductivity equation; anti-plane elasticity; stress; gradient blow-up; extreme conductivity; stress concentration factor}

%%%%%%%%%%%%%%%%%%%%%%%%%%%%%%%%%%%%%%%%%%%%%%%%%%%%%%%%%%%%%%%
\section{Introduction and statement of results}
%%%%%%%%%%%%%%%%%%%%%%%%%%%%%%%%%%%%%%%%%%%%%%%%%%%%%%%%%%%%%%%

In presence of closely located stiff inclusions embedded in the relatively weak matrix, high stress concentration occurs in the narrow region between two inclusions. Such a phenomenon typically occurs in fiber-reinforced materials and the stiff inclusions represent the cross section of fibers. Recently, much effort has been devoted to quantitative understanding of this stress concentration. In this paper we continue our investigation on this and establish an efficient method to compute the magnitude of the stress concentration that immediately yields an asymptotic formula for the stress distribution and an optimal estimate for the concentration.

To describe the problem and results in a precise manner, let $D_1^0$ and $D_2^0$ be a pair of (touching) bounded domains with $\mathcal{C}^{2, \Gg}$ ($\Gg>0$) boundaries such that
\beq
D_1^0 \subset \{ (x,y) \in \Rbb^2 ~|~ x <0\}, \quad D_2^0 \subset \{ (x,y) \in \Rbb^2 ~|~ x >0\},
\eeq
\beq
\p D_1^0 \cap \p D_2^0 = \{(0,0)\},
\eeq
and $D_1^0$ and $D_2^0$ are convex at $(0,0)$. The domains $D_1^0$ and $D_2^0$ are strongly convex at $(0,0)$ if both $D_1^0$ and $D_2^0$ have positive curvatures there. By translating $D_2^0$ by a positive number $\Ge$ along $x$-axis, while $D_1^0$ is fixed,  we obtain $D_2^\Ge$, {\it i.e.},
\beq
D_2^\Ge := D_2^0 + \left(\Ge , 0\right).
\eeq
When there is no possibility of confusion, we drop superscripts and denote
\beq
D_1:=D_1^0, \quad D_2:=D_2^\Ge.
\eeq

For a given harmonic function $h$ in $\mathbb{R}^2$, let $u_\Ge$ be the solution to the problem
\beq \label{gov_eq}
\begin{cases}
\ds \Delta u_\Ge = 0 \quad&\mbox{in } \Rbb^2 \setminus \overline{(D_{1} \cup D_{2})},\\
\ds u_\Ge = \Gl_{j} \ (\mbox{constant}) \quad&\mbox{on } \p D_{j},~j=1,2, \\
\ds u_\Ge(\Bx) - h(\Bx) =O(|\Bx|^{-1}) ~ &\mbox{as }
|\Bx|\rightarrow \infty,
\end{cases}\eeq where the constants $\Gl_{i }$ are determined by the
conditions
\beq\label{intzero}
 \int_{\p D_{1}} \p_\nu u_\Ge ds = \int_{\p D_{2}} \p_\nu u_\Ge ds =0 .
\eeq
Here and throughout this paper $\p_\nu u_\Ge$ denotes the outward normal derivative of $u_\Ge$ on $\p D_{j}$ $(j=1,2)$. It is worth emphasizing that the constants $\Gl_1$ and $\Gl_2$ may or may not be different depending on the given $h$.

As mentioned before, inclusions $D_1$ and $D_2$ represent the two dimensional cross-sections of two parallel elastic fibers embedded in an infinite elastic matrix and $\Ge$ is the distance between them.
The solution $u_\Ge$ represents the out-of-plane elastic displacement, and $\nabla u_\Ge$ is proportional to the shear stress. The problem \eqnref{gov_eq} may also be regarded as two dimensional conductivity equation in which case $D_1$ and $D_2$ represent perfect conductors of infinite conductivity.
It is worth mentioning that we consider the situation where there are only two inclusions since our interest lies in estimating local high concentration of stress in the narrow region between two inclusions. There is a study to estimate global stress in a composite (with many inclusions) using a network approximation. We refer to \cite{BKN} and references therein for that. We also mention that the problem under consideration in this paper has some connection with effective properties of composites with highly conducting inclusion. See \cite[Section 10.10]{milton_book} for this connection.

In general, $\nabla u_\Ge$ becomes arbitrarily large as the distance $\Ge$ between two inclusions tends to zero, and the problem is to derive pointwise estimates of $\nabla u_\Ge$ in terms of $\Ge$. This problem was raised in \cite{bab} and there has been significant progress on it. It has been proved that the generic blow-up rate of $\nabla u_\Ge$ is $1/\sqrt{\Ge}$ in two dimensions \cite{AKLLL, AKLLZ, AKLim, BLY, BC, keller, LY2, Y, Y2}, and $|\Ge \ln \Ge|^{-1}$ in three dimensions \cite{BLY, BLY2, KLY13, lekner10, lekner11, lekner, LY}. We emphasize that the gradient may or may not blow up depending on the given background harmonic function $h$. For example, in the configuration of this paper the gradient blows up if $h(x,y)=x$ and it does not if $h(x,y)=y$ for circular inclusions $D_1$ and $D_2$. It is worth while to mention that the insulating case in two dimensions can be treated by duality as done in \cite{AKLim} for example. But the insulating case in three dimensions is an open problem: it is not even clear if the gradient actually blows up in three dimensions. It is also worth while to mention that if the conductivity of the inclusions is finite (away from $\infty$ and $0$), $\nabla u_\Ge$ is bounded regardless of $\Ge$ \cite{bv, ln, lv}.

Recently a better understanding of the stress concentration has been obtained: an asymptotic behavior of $\nabla u_\Ge$ has been characterized by the singular function associated with $D_1$ and $D_2$, as $\Ge$ tends to $0$. The singular function, denoted by $q_\Ge$, is the solution to the following problem:
\beq \label{q_equ}
\begin{cases}
\ds \Delta q_\Ge = 0 \quad&\mbox{in } \Rbb^2 \setminus \overline{(D_{1} \cup D_{2})},\\
\ds q_\Ge = \mbox{constant} \quad&\mbox{on } \p D_j, ~ j=1,2, \\
\ds q_\Ge(\Bx)  =O(|\Bx|^{-1}) ~ &\mbox{as }
|\Bx|\rightarrow \infty,\\
\ds \int_{\p D_{1}} \p_\nu q_\Ge ds =-1, \quad  &\ds\int_{\p D_{2}} \p_\nu q_\Ge ds =1 .
\end{cases}
\eeq
We emphasize that the constant values of $q_\Ge$ on $\p D_1$ and $\p D_2$ are different, so that $\nabla q_\Ge$ blows up as $\Ge \to 0$.

Let us recall some important facts about $q_\Ge$: If $D_1^0$ and $D_2^0$ are disks, then $q_\Ge$ is given explicitly  by
 \beq\label{singular}
 q_\Ge(\Bx):= \frac{1}{2\pi} \left( \ln |\Bx-\Bp_1|- \ln |\Bx-\Bp_2|\right),
 \eeq
where $\Bp_1 \in D_1$ is the fixed point of the mixed reflection $R_1 R_2$ where $R_j$
is the reflection with respect to $\p D_j$, $j=1,2$, and $\Bp_2 \in D_2$ is that of $R_2 R_1$. We emphasize that these points can be computed easily (see \eqref{p_i_esti}). More generally, if $D_1^0$ and $D_2^0$ are strongly convex at $(0,0)$, then let $B_1$ and $B_2$ be disks osculating to $D_1$ and $D_2$ at $(0,0)$ and $(\Ge, 0)$, respectively, and let $q_{B,\Ge}$ be the singular function associated with $B_1$ and $B_2$ as given in \eqnref{singular}. Then, it is proved in \cite{ACKLY} that the behavior of $\nabla q_{\epsilon}$ is almost explicitly  described as
\beq\label{qBGe}
\nabla q_\Ge = \nabla q_{B,\Ge}(1+ O(\Ge^{\Gg/2})) + O(1)
\eeq
when $\p D_j$ is $\mathcal{C}^{2,\Gg}$ and $\gamma \in (0,1)$.
Thus it follows that
\beq
\norm{\nabla q_\Ge}_{\infty} = \frac {\sqrt {\kappa_1 + \kappa_2} } { \sqrt 2 \pi }  \frac 1 {\sqrt{\Ge}} + O(1),
\eeq where $\kappa_2$ and $\kappa_2$ are the curvatures of $D_1 ^0$ and $D_2 ^0$ at $(0,0)$, respectively.

Using the singular function $q_\Ge$, the solution $u_\Ge$ to \eqnref{gov_eq} can be decomposed as
\beq\label{decom1}
u_\Ge = \Ga_{\Ge} q_\Ge + r_\Ge
\eeq
where
\beq
\Ga_{\Ge} = \frac {u_\Ge |_{\p D_2}  - u_\Ge |_{\p D_1}}  { q_\Ge |_{\p D_2} - q_\Ge |_{\p D_1}}.
\eeq
Observe that $r_\Ge$ is also constant on $\p D_1$ and $\p D_2$, and $r_\Ge |_{\p D_1} = r_\Ge |_{\p D_2}$, so that  $\nabla r_\Ge$ is bounded on bounded subsets of $\Rbb^2 \setminus \overline{(D_{1} \cup D_{2})}$ (see \cite{KLY}).
It means that the term $\Ga_{\Ge} \nabla q_\Ge$ is responsible for the blow-up of $\nabla u_\Ge$, or more precisely,
\beq\label{uGeqe}
\nabla u_\Ge = \Ga_{\Ge} \nabla q_\Ge + O(1) \quad\mbox{as } \Ge \rightarrow 0.
\eeq
In particular, $\Ga_\Ge$ represents the magnitude (normalized by $|\nabla q_\Ge|$) of the blow-up. So, it is appropriate to call the constant $\Ga_\Ge$ the {\it stress concentration factor}.

The purpose of this paper is to analyze the stress concentration factor $\Ga_{\Ge}$ when $D_1^0$ and $D_2^0$ are convex at $(0,0)$. We are particularly interested in finding $\lim_{\Ge \to 0} \Ga_\Ge$ (existence of the limit is a part of the study).

There have been some work on the stress concentration factor. It is proved in \cite{KLY} that if $D_1^0$ and $D_2^0$ are disks, then
\beq\label{concent_fac_disk}
\Ga_\Ge = \frac{2 r_1 r_2}{r_1 + r_2} (\Bn \cdot \nabla h) (0,0) + O(\sqrt{\Ge}) \quad\mbox{as } \Ge \to 0,
\eeq
where $r_j$ is the radius of $D_j$, $j=1,2$, and $\Bn$ is the outward unit normal vector to $\p D_1$ at $(0,0)$. An estimate for $\Ga_\Ge$  in terms of curvatures, size and $\Ge$  was  established in  \cite{LY2} under the assumption that an inclusion has a much higher curvature than its size.  It is also proved in \cite{ACKLY} that if $D_1^0$ and $D_2^0$ are strongly convex at $(0,0)$, then
\beq\label{hqGe}
\Ga_{\Ge} = \frac{\sqrt 2 \pi} {\sqrt {\Gk_1 + \Gk_2}} \frac 1 {\sqrt \Ge } \int_{\p D_1 \cup \p D_2} h \p_{\nu} q_{\Ge} ds \left( 1 + O(\Ge^{\Gg/2}) \right),
\eeq
and, as a consequence, that $\Ga_\Ge$ is bounded regardless of $\Ge$.

Observe that even if \eqnref{qBGe} yields a good information of $q_\Ge$ on the narrow region in between two inclusions, it is still difficult to evaluate the integral on the righthand side of \eqnref{hqGe} since it requires global information of $\p_\nu q_\Ge$.
In this paper, we present a new efficient method for finding $\lim_{\Ge \to 0} \Ga_\Ge$. It turns out that the limit is given as a certain integral of the solution $u_0$ for the touching case, namely,
\beq\label{u_0}
\begin{cases}
\ds \Delta u_{0}  = 0 \quad&\mbox{in } \GO,\\
\ds u_{0} = \Gl_0 \quad&\mbox{on } \p \GO, \\
\ds u_{0}(\Bx) - h(\Bx) =O(|\Bx|^{-1}) \quad&\mbox{as } |\Bx| \to \infty,
\end{cases}
\eeq
where $\GO:= \Rbb^2 \setminus \overline{(D_1^{0} \cup D_2^{0})}$ and $\Gl_0$ is a constant determined by the additional condition
\beq
\int_{\Omega} |\nabla (u_0 - h)|^2 d A< \infty.
\eeq

To obtain the main result of this paper (Theorem \ref{1st_main}), we first consider the touching case problem \eqnref{u_0}.
We show that there exists a unique solution $u_0$ to \eqnref{u_0} and $\nabla u_0(x,y)$ decays exponentially fast as $(x,y)$ approaches to $(0,0)$ in $\GO$ (Theorem \ref{thmtouching}). It is worth mentioning that $\GO$ has cusps at the origin. We also prove a similar theorem for the residual part $r_\Ge$ in \eqnref{decom1} (Theorem \ref{asymdecay2}). This result was also obtained in \cite{LLBY} in a more general context. However, we include a proof in this paper since it is completely different from that in the paper mentioned above.

We prove these results in somewhat more general setting: We assume that the domains $D_1^0$ and $D_2^0$ are convex at $(0,0)$ and their order of contact at the point is $2m$ for some positive integer $m$. Thus, if $\p D_j^{0}$ near $(0,0)$ is given as the graph of $x=x_j(y)$ ($j=1,2$), then there are constants $\Gd_0>0$ and  $c_j>0$, $j=1,\ldots,4$,  such that
\beq\label{contact}
-c_1 y^{2m} \le x_1(y) \le -c_2 y^{2m} \quad\mbox{and}\quad c_3 y^{2m}\le x_2(y) \le c_4 y^{2m},
\eeq
for $|y| < \Gd_0$. If $D_1^0$ and $D_2^0$ are strongly convex at $(0,0)$, then $m=1$.

In terms of the solution $u_0$ to \eqnref{u_0} we obtain the following theorem regarding the asymptotic behavior of the stress concentration factor.
\begin{thm}\label{1st_main}
Suppose that $\p D_j^{0}$ $(j=1,2)$ are $\mathcal{C}^{2,\Gg}$ for some $\Gg>0$, $D_1^0$ and $D_2^0$ are convex at $(0,0)$, and their order of contact at $(0,0)$ is $2m$. Let $u_0$ be the solution to \eqref{u_0} and let
\beq
\Ga_0 := \int_{\p D_{1}^0} \p_{\nu} u_{0} \, ds .
\eeq
Then,
\beq\label{eq:limit}
\Ga_\Ge = \Ga_0 + O \left( \Ge|\log \Ge |^{2m-1}\right)~~\mbox{ as }\Ge \rightarrow 0.
\eeq
\end{thm}

As an immediate consequence of Theorem \ref{1st_main} and \eqnref{uGeqe}, we obtain
\beq\label{mainsing2}
\nabla u_\Ge = \Ga_{0} \nabla q_\Ge + O(1) \quad\mbox{as } \Ge \rightarrow 0
\eeq
in any bounded subset of $\mathbb{R}^2 \setminus (D_1 \cup D_2)$. Thus, the limit $\Ga_0$ can be regarded as an alternative concentration factor. Moreover, if $D_1^0$ and $D_2^0$ are strongly convex at $(0,0)$, then we have from \eqnref{singular} and \eqnref{qBGe} that
\beq\label{mainsing}
\nabla u_\Ge(\Bx) = \frac{\Ga_0}{2 \pi} \left( \frac {\Bx - \Bp_1} {|\Bx - \Bp_1|^2} -   \frac {\Bx - \Bp_2} {|\Bx - \Bp_2|^2} \right)(1+ O(\Ge^{\Gg/2})) + O(1).
\eeq

These formulas have some important consequences. As a first consequence, we have the following identity (see section 4 for a proof):
\beq\label{limit}
\lim_{\Ge \to 0} \sqrt{\Ge} |\nabla u_\Ge(\Ge/2, 0)| = \frac {\Ga_0\sqrt {\Gk_1 + \Gk_2 }}{\sqrt 2 \pi},
\eeq where $\Gk_i$ is the curvature of $\p D_i ^0$ at $(0,0)$ for $i=1,2$. Note that $(\Ge/2, 0)$ is a point where $|\nabla u_\Ge|$ has a value close to the maximal concentration.

Another consequence of \eqnref{mainsing2} and \eqnref{mainsing} is related to numerical computation of $\nabla u_\Ge$. Since high concentration of the gradient occurs in the narrow region, fine meshes may be required to compute $\nabla u_{\Ge}$. However, since \eqnref{mainsing2} and \eqnref{mainsing} extract the major singular term in an explicit way, it suffices to compute the residual term $\nabla r_\Ge$ for which only regular meshes are required. This idea was exploited in \cite{KLY} in the special case when $D_j$'s are disks using \eqref{singular},  \eqref{decom1} and \eqref{concent_fac_disk}. Implementation of this idea for the general case of strongly convex domains will be the subject of the forthcoming work. It is worth mentioning that there are some other methods to compute the solution when $D_1$ and $D_2$ are disks. See, for examples, \cite{CG, MPM}.

The last subject of this paper is regarding computation of $\Ga_0$. It turns out that, thanks to exponentially decaying property of the solution to the touching case, $\Ga_0$ can be computed numerically only using regular meshes by truncating the narrow region near $(0,0)$.

This paper is organized as follows. We investigate the touching case in section \ref{cusp}. In section \ref{sec3} we obtain an estimate for the gradient of the residual term $r_\Ge$. Section \ref{sec:proof} is to prove Theorem \ref{1st_main}. In the last section we present a way to compute good approximations of $\Ga_0$.

%%%%%%%%%%%%%%%%%%%%%%%%%%%%%%%%%%%%%%%%%%%%%%%%%%%%%%%%%%%%%%%
\section{The solution for the touching case}\label{cusp}
%%%%%%%%%%%%%%%%%%%%%%%%%%%%%%%%%%%%%%%%%%%%%%%%%%%%%%%%%%%%%%%

In this section we prove the following theorem regarding the problem \eqnref{u_0}.
\begin{thm}\label{thmtouching}
Suppose that $\p D_j^{0}$ $(j=1,2)$ are $\mathcal{C}^{2,\Gg}$ for some $\Gg>0$, $D_1^0$ and $D_2^0$ are convex at $(0,0)$, and their order of contact at $(0,0)$ is finite. Then, there is a unique solution $u_0$ to \eqref{u_0}, and there are positive constants $A$, $C$ and $\Gd$ such that
\beq\label{asymdecay}
|\nabla u_0(x,y)| \le C \exp \left( - \frac{A}{|y|} \right)
\eeq
for $|y| \leq \Gd$ and $x_1(y) < x < x_2 (y)$, where $x_1$ and  $x_2$ are the defining functions of $\p D_1^0$ and $\p D_2^0$ near $(0,0)$.
\end{thm}

The estimate \eqnref{asymdecay} follows from
\beq\label{u0decay}
|u_0(x,y)- \Gl_0| \le C \exp \left( - \frac{A}{|y|} \right)
\eeq
by a standard estimate for harmonic functions. Here the constants $A$ and $C$ may differ at each occurrence. In fact, since $\p D_j^{0}$ are $\mathcal{C}^{2,\Gg}$ and $u_0-\Gl_0=0$ on $\p D_j^0$, one can show that $u_0(x,y)-\Gl_0$ can be extended by reflection (after the conformal transformations to outside a disk) as harmonic functions for $(x,y)$ satisfying
$$
x_1(y) - s y^{2m} < x < x_2 (y) $$
and
$$
x_1(y) < x < x_2 (y) + s y^{2m}
$$ for some $s>0$, and the same estimate \eqnref{u0decay} holds for the extended functions. Here $2m$ is the order of contact at $(0,0)$. So, for each $(x,y) \in \Rbb^2 \setminus (D_1^0 \cup D_2^0)$, there is $r>0$ such that
$r > t y^{2m}$ for some $t>0$ and $u_0$ is harmonic in $\overline{B_r(x,y)}$. ($B_r(c)$ denotes the disk of radius $r$ with the center at $c$.) So, we have
\begin{equation}
|\nabla u_0(x,y)| \le \frac{C}{r^3} \int_{B_r(x,y)} |u_0-\Gl_0|~ dA \le \frac{C}{y^{6m}} \exp \left( - \frac{A}{|y|} \right) \le C' \exp \left( - \frac{A'}{|y|} \right) \label{standa_esti}
\end{equation}
for some constant $A'$ and $C'$, which is the desired estimate.

The rest of this section is devoted to proving existence and uniqueness of $u_0$, and \eqnref{u0decay}.
To construct the solution to \eqnref{u_0} we use the transformation $1/z$, following \cite{MS}.
Identify $\Bx=(x,y)$ in the plane with $z=x+iy$ and let
$$
\Phi(z)= \frac{1}{z}.
$$
Define
$$
\widetilde{\GO}:= \Phi(\GO), \quad \GG_1:= \Phi(\p D_1^0), \quad \GG_2:= \Phi(\p D_2^0).
$$
Note that $\GG_1$ and $\GG_2$ are simple curves lying in the left and right half spaces, respectively, and $\widetilde{\GO}$ is the region enclosed by $\GG_1$ and $\GG_2$. Since $\p D_j^0$ is $\mathcal{C}^{2,\Gg}  (\Gg>0)$ and $D_j^0$ is convex at $(0,0)$ for $j=1,2$, one can easily see that there are constant $a<b$ such that
\beq\label{GOsubset}
\widetilde{\GO} \subset \{ w=\xi+i \eta ~|~ a < \xi < b \}.
\eeq
Moreover, $\GG_1$ near $\infty$ is given by $\xi=\psi_1(\eta)$ for some function $\psi_1$ satisfying
\beq\label{psione}
\psi_1(\eta) \le  -C_1 |\eta|^{2-2m}
\eeq
for some constant $C_1>0$, and $\GG_2$ near $\infty$ is given by $\xi=\psi_2(\eta)$ for some function $\psi_2$ satisfying
\beq\label{psitwo}
\psi_2(\eta) \ge  C_2 |\eta|^{2-2m}
\eeq
for some constant $C_2>0$. In fact, we have
\beq\label{psioneeta}
\psi_1(\eta) = \frac{x_1 (y)}{y^2 + x_1(y)^2} \quad\mbox{with} \quad \eta=\frac{-y}{y^2 + x_1(y)^2}
\eeq
on $\p D_1^0$ near $(0,0)$.
Thanks to \eqnref{contact}, we have
$$
a<\psi_1(\eta) \le  \frac{-c_1 y^{2m}}{y^2 + x_1(y)^2} \le -C_1 |\eta|^{2-2m}.
$$
Thus we have \eqnref{psione}. \eqnref{psitwo} can be proved similarly.

We need the following lemma whose proof will be given after completing the proof of Theorem \ref{thmtouching}.
\begin{lem}\label{lem:decay1}
Let $\psi_j$ $(j=1,2)$ be as defined by \eqref{psioneeta}, and let $a$ and $b$ be the constants such that
\beq
a  < \psi_1(\eta) < \psi_2(\eta) < b
\eeq
for all $\eta >L $, where $L$ is a large number. Let $R$ be a domain given by
\beq
R:= \{ (\xi, \eta) ~|~ \eta >L, \ \psi_1(\eta) < \xi < \psi_2(\eta) \},
\eeq
and let $U$ be the solution in $H^1(R)$ to the problem
\beq \label{U}
\left\{
\begin{array}{rll}
\Delta U &=0 \quad &\mbox{in } R, \\
U&=0\quad &\mbox{on } \xi=\psi_j(\eta), \ j=1,2, \\
U&=\Gvf \quad &\mbox{on } \GG:=\{ (\xi,L) ~|~ \psi_1(L) <  \xi < \psi_2(L) \},
\end{array}
\right.
\eeq
where $\Gvf$ is a bounded function. Then there are positive constants $A$ and $C$ such that
\beq\label{Uxieta}
|U(\xi, \eta)| \le C e^{-A\eta}
\eeq
for all  $\eta >L$.
\end{lem}

Because of \eqnref{GOsubset}, the Poincar\'e inequality holds in $\widetilde{\GO}$: for all $\tilde u \in H^1_0(\widetilde{\GO})$ (the standard Sobolev space with the zero trace)
\begin{align} \label{poincare}
\norm{\tilde u}^2 _{L^2 ( \widetilde \Omega)} \le C \norm{\nabla \tilde u}^2 _{L^2 ( \widetilde \Omega)}
\end{align}
for some constant $C$. So, one can apply the Lax-Milgram Theorem to show that for $f \in H ^{-1}(\widetilde \Omega)$ there exists a unique solution $\tilde v \in H^1_0( \widetilde{\GO} )$ to
\beq \label{tilde_v}
\left\{
\begin{array}{rl}
\ds \Delta \tilde v   &= f \quad\mbox{in } \widetilde \Omega,\\
\ds \tilde v  &= 0 \quad\mbox{on } \p \widetilde \Omega.
\end{array}
\right.
\eeq

We choose $r_0>0$ such that
$$
D _{1}^0\cup D_{2}^0\subset B_{r_0/2}(0) .
$$
Let $\chi$ be a smooth function such that
$\chi (z) = 1$ if $z \in B_{r_0} (0)$ and $\chi(z) = 0$ if $z \notin B_{2r_0}(0)$.
For $h$ given in \eqnref{u_0}, let
$$
f (w)= \Delta_w \left(\chi \left(\frac 1 w \right)   h \left(\frac 1 w \right)\right) ,
$$
and let $\tilde v$ be the solution to \eqnref{tilde_v} with this $f$. Then one can check that $u_0$ given by
\beq
u_{0}(z) = h(z)+\left({\tilde v}\left( \frac 1 z \right)  -\chi(z)h (z) - {\tilde v} ( 0 ) \right)
\eeq
is the solution to \eqnref{u_0} and the constant value $\Gl_0$ is given by $ - {\tilde v} ( 0 )$. The uniqueness of the solution follows easily from the maximum principle.

Now, we show \eqref{u0decay}. If $z \in B_{r_0}(0) \setminus \overline{D _{1}^0\cup D_{2}^0}$, then we have
\beq
u_0(z) -\Gl_0 = {\tilde v}\left( \frac 1 z \right). \label{u_0v}
\eeq

Choose $L$ so large that the support of $f$ lies in between two lines $\eta = \pm L$. Let $\widetilde\GO_{\pm L} := \widetilde\GO \cap \{ \pm \eta >L \}$, respectively.  The boundedness of $ \tilde v (\xi \pm (L+1)i)$ can be shown easily by a standard estimate for harmonic functions similarly to \eqref{standa_esti}, since $\tilde v  = 0 $ on $\p \widetilde\GO_{\pm L} \cap \p \widetilde \Omega $ and $\tilde v \in L^2 (\widetilde{\GO})$.  We thus apply Lemma \ref{lem:decay1}  to obtain
\beq\label{vdecay}
|\tilde v(\xi+i\eta)|  \le C_1 e^{-A_1 |\eta|} \quad\mbox{for}~ |\eta| > L+1
\eeq
for some positive constant $A_1$ and $C_1$. We may choose  a small positive number $\Gd_1$ so that
$\Phi(x,y) \in \widetilde\GO_{+ (L+1)} \cup  \widetilde\GO_{- (L+1)}$ for all $x+iy$ satisfying
$|y| < \Gd_1$ and $x_1(y) < x < x_2 (y)$. Then, by \eqref{u_0v} and \eqref{vdecay}, we have
\begin{align}  \label{esti_u_0}
&\left| u_0(x,y) + {\tilde v} ( 0 )\right|= \left|  {\tilde v} \left(\xi+\eta i \right) \right| \le C_2 e^{-A_1|\eta|}  \le C_3 e^{-\frac{A_2}{|y|}},
\end{align}
for $ |y| < \Gd_1$. The last inequality follows from \eqnref{psioneeta}, since $|x_1(y)| \simeq |y|^{2m}$. Here and throughout this paper, $a \simeq b$ stands for $\frac 1 C a \leq b \leq C a$ for some constant $C$ independent of $\Ge$.
 This completes the proof of Theorem \ref{thmtouching}.
\qed

\medskip

\noindent{\sl Proof of Lemma \ref{lem:decay1}}.
By translating and scaling if necessary, we may assume $a=0$, $b=\pi$  and $L=0$. Let
$$
\tilde R:= \{ (\xi, \eta) ~|~ \eta >0, \ 0 < \xi < \pi \}.
$$
Decompose $\Gvf$ as $\Gvf=\Gvf_+ - \Gvf_-$ where $\Gvf_\pm $ are nonnegative and bounded, and then extend $\Gvf_\pm$ to $[0, \pi]\times \{0\}$ by assigning $0$ outside $\GG$, and denote them by $\tilde \Gvf_\pm$. Let $V_\pm$ be a solution in $H^1 (\tilde R)$ to
\beq
\left\{
\begin{array}{ll}
\ds \Delta V_{\pm}   = 0 \quad&\mbox{in } \tilde R,\\
\ds V_{\pm} (0, \eta)= V_{\pm} (\pi, \eta) = 0, \quad & \eta >0, \\
\ds V_{\pm}(\xi, 0)   =\tilde \Gvf_\pm(\xi),  \quad & 0 \le \xi \le \pi.
\end{array}
\right. \notag
\eeq
Since $\tilde \Gvf_\pm \ge 0$, we have $V_\pm \ge 0$, and by the maximum principle, we have
\beq\label{V-UV+}
-V_{-} \leq U  \leq V_{+} \quad \mbox{in } \tilde R.
\eeq
One can find the solutions $V_{\pm}$ by separation of variables. In fact, we have
$$
V_{\pm} (\xi,\eta) = \sum_{n=1}^{\infty} a_n^\pm \sin n \xi \, e^{-n\eta},
$$
where $a_n^\pm$ is the Fourier coefficients of $\Gvf_\pm$. In particular, we have
\beq
|V_{\pm} (\xi,\eta)| \leq \left(\sum_{n=1}^{\infty} |a_n^\pm|^2 \right)^{1/2} \left(\sum_{n=1}^{\infty} e^{-2n\eta} \right)^{1/2} \le C e^{-\eta}. \label{Vplus}
\eeq for $\eta \geq 1$. Even for  $0< \eta < 1$, this inequality holds with another constant $C$ since $\tilde\varphi_\pm$ are bounded. Thus, \eqnref{Uxieta} follows from \eqnref{V-UV+}.
This completes the proof.
\qed

%%%%%%%%%%%%%%%%%%%%%%%%%%%%%%%%%%%%%%%%%%%%%%%%%%%%%%%%%%
\section{The behavior of $\nabla r_\Ge$ in the narrow region} \label{sec3}
%%%%%%%%%%%%%%%%%%%%%%%%%%%%%%%%%%%%%%%%%%%%%%%%%%%%%%%%%%
In this section, we consider the behavior of the gradient of $r_\Ge$ given in \eqref {decom1} in the narrow region between $D_1$ and $D_2$ which we denote by $N_\Gd$ for $\delta >0$, namely,
\beq
N_\Gd:= \left\{(x,y) ~|~ x_1 (y)<x< x_2 (y) +\Ge, \ |y|< \Gd  \right\}.
\eeq

Recall that $r_\Ge$ satisfies
\beq
\begin{cases}
\ds \Delta r_\Ge  = 0  \quad\mbox{in } \GO:= \mathbb{R}^2 \setminus \overline{ D_1 \cup D_2} ,\\
\ds r_\Ge |_{\p D_1} =r_\Ge |_{\p D_2} = \mbox{constant} ,\\
\ds r_\Ge (\Bx)- h (\Bx) = O (|\Bx|^{-1})  \quad\mbox{as } |\Bx| \rightarrow \infty.
\end{cases}
\eeq
In the previous section, it has been shown that $\nabla u_0$ is decreasing exponentially near origin. The following theorem shows that $\nabla r_\Ge$ has such a decay property. As mentioned in Introduction, this result was also obtained in \cite{LLBY} in a more general setting. But two proofs are completely different. 

\begin{thm}\label{asymdecay2}
Suppose that $\Ge$ is sufficiently small. There are positive constants $A$, $C$ and $\Gd$ independent of $\Ge >0$ such that
\beq\label{mainest}
| \nabla r_{\Ge} (x,y) | \leq C \exp \left(-  \frac A {\sqrt {\Ge} + |y|}\right)
\eeq
for any $(x,y) \in N_\Gd$.
\end{thm}

We prove the following lemma from which Theorem \ref{asymdecay2} follows by a standard elliptic estimate as explained briefly in the previous section.
\begin{lem}\label{d**}
Suppose that $\Ge$ is sufficiently small. There are positive constants $A$, $C$ and $\Gd$ independent of $\Ge$ such that
\beq\label{mainest2}
| r_{\Ge} (x,y) - r_\Ge |_{\p D_1} | \leq C \exp \left(-  \frac {A}{\sqrt {\Ge} + |y|}\right)
\eeq
for any $(x,y) \in N_\Gd$.
\end{lem}
\pf
As before we identify points $(x,y)$ in $\Rbb^2$ with $z= x+iy$ in $\Cbb$. Choose two disks $B_1$ and $B_2$ whose centers are on the real axis such that
$$
B_j \subset D_j ,~j=1,2, ~~ \p B_1 \cap \p D_1 = \{0\}, ~~\mbox{and} ~~\p B_2 \cap \p D_2 = \{\Ge\}.
$$
Let $c_j$ and $\rho_j$ be the center and radius, respectively, of $B_j$ for $j=1,2$. It is convenient to assume that $c_2=1+\Ge$ so that $\rho_2=1$.  Let
\beq
\Phi_1 ( z) = \frac 1 { z - ( 1 + \Ge) }
\eeq
which is the reflection with respect to $\p B_2$ (and translation), and
\beq
\GO_1 : =  \Phi_1  \left(\GO  \right), \quad B_3:=\Phi_1(B_1), \quad B_4=\Phi_1(B_2),
\eeq
and let $c_j$ and $\rho_j$ be the center and radius of $B_j$. Then, $c_4=0$, $\Gr_4=1$, and $\Gr_3= c_3-\Phi_1(0)$. Observe that
\beq\label{domain_subset}
\GO_1 \subset B_4 \setminus  \overline{ B_3},
\eeq
the reflected domain $\GO_1$ touches $\p B_3$ and $\p B_4$ at $\Phi_1(0)$ and $\Phi_1(\Ge)$, respectively, and
\beq
\mbox{dist} \, (\p B_3, \p B_4 ) = \Phi_1(0)-\Phi_1(\Ge) = - \frac{1}{1+\Ge} +1 = \Ge+ O(\Ge^2). \label{distB_3B_4}
\eeq

Let
\beq\label{setS}
S : = \left\{w ~|~ w =\xi+i\eta \in B_4 \setminus  \overline{B_3}, ~ \xi <c_3, ~ |\eta| < \Gr_3/2  \right\}.
\eeq
Then one can choose $\Gd$  independently of $\Ge$ so that $\Phi_1(N_\Gd) \subset S$. Let
\beq
\tilde r_\Ge(w) := r_\Ge \circ \Phi_1^{-1}(w) - r_\Ge|_{\p D_1}, \quad w \in \GO_1. \label{tilderGe2}
\eeq
If $z=x+iy \in N_\Gd$ and $w=\xi+i\eta = \Phi_1(z)$, then $\eta \simeq y$. Thus in order to prove \eqnref{mainest2}, it suffices to show
\beq\label{mainest3}
| \tilde r_\Ge(w) | \leq C \exp \left(-  \frac {A}{\sqrt {\Ge} + |\eta|}\right)
\eeq
for any $w \in \Phi_1(N_\Gd)$.

We now transform $B_3$ so that the transformed disk becomes concentric to $B_4$ ($B_4$ is the unit disc). For that purpose let us write a lemma which can be easily verified.
\begin{lem}
Let $B_\Gr(c)$ be a disk such that $\overline{B_\Gr(c)} \subset B_1(0)$. Then there is $\Ga$ with $|\Ga| <1$ and $\Gr_*>0$ such that the M\"obius transform $\Gvf_\Ga$ defined by
\beq\label{Phi}
\Gvf_\Ga(w) = \frac{w-\Ga}{1-\bar{\Ga} w}
\eeq
maps $B_\Gr(c)$ onto $B_{\Gr_*}(0)$. In fact, $\Ga$ is given by
\beq\label{alpha_def}
\Ga= \left[ (|c|^2 -\Gr^2 +1)- \sqrt{(|c|^2 -\Gr^2 +1)^2- 4|c|^2} \right] \frac{c}{2|c|^2}.
\eeq
\end{lem}

It is worth mentioning that M\"obius transforms are automorphisms on $B_{1}(0)$.

\par  Let $\Phi_2$ be the M\"obius transform defined by \eqnref{Phi} and \eqnref{alpha_def} with $c=c_3$ and $\Gr= \Gr_3$, and let $B_5=B_{\Gr_5}(c_5):= \Phi_2(B_3)$. Then $c_5=0$. Since $\Gr_3=c_3+ \frac{1}{1+\Ge}$, one can see from \eqnref{alpha_def} that $\Ga$ is real and satisfies
\beq\label{Ga=-1}
\Ga = -1 + \Gb\sqrt{\Ge} + (1+ \frac{1}{c}) \Ge + O(\Ge \sqrt{\Ge}),
\eeq
where
$$
\Gb=  \sqrt{\frac{2 (c_3+1)}{|c_3|}}.
$$
To compute $\Gr_5$, we observe that $\Phi_2(-\frac{1}{1+\Ge}) \in \p B_5$, and from \eqnref{Ga=-1} that
$$
\Phi_2(-\frac{1}{1+\Ge}) = \frac{-\frac{1}{1+\Ge} - \Ga}{1 +\frac{\overline{\Ga}}{1+\Ge}} = -1 + \Gg \sqrt{\Ge} + O(\Ge),
$$
where $\Gg= 2/\Gb$. So, we have
\beq\label{Gr_5=1}
\Gr_5= 1-\Gg \sqrt{\Ge} + O(\Ge).
\eeq
We emphasize that \eqnref{Gr_5=1} implies in particular that
\beq
\mbox{dist} \, (\p B_5, \p B_4 ) = \Gg \sqrt{\Ge} + O(\Ge),
\eeq
since $ B_4  =  B_1 (0)$.

\par The proof of the following lemma will be given later in this section. Here $\mbox{arg} (z)$ for $z\neq 0$ is supposed to take a value in $[0,2\pi)$.
\begin{lem}\label{C_*}
Suppose that $\Ge$ is sufficiently small. There exists a constant $C>0$ independent of $\Ge$ such that
\beq\label{C*eq}
\mbox{arg} \left( \Phi_2 \left(w \right) \right) \geq   \frac {C\sqrt \Ge}{|\eta|+ \sqrt \Ge}
\eeq
for $w  = \xi + \eta i \in S$ with $\eta \ge 0$, and
\beq
2\pi -\mbox{arg} \left( \Phi_2 \left(w \right) \right) \geq   \frac {C\sqrt \Ge}{|\eta|+ \sqrt \Ge}
\eeq
for $\eta \le 0$.
\end{lem}

Let us introduce one more transformation $\Phi_3$:
\beq
\Phi_3(\Gz)= \log \Gz
\eeq
with the branch cut on the positive real axis. Then $\Phi_3$ maps $(B_4\setminus \overline{B_5}) \setminus \{\mbox{positive real axis} \}$ onto the rectangle $(a_0, 0) \times (0, 2\pi)$ where $a_0= \log \Gr_5 <0$. We emphasize that
\beq\label{a=-Gg}
a_0 =- \Gg\sqrt{\Ge} + O(\Ge),
\eeq
which is a consequence of \eqnref{Gr_5=1}.

Let $\Gt_0$ be the constant on the righthand side of \eqnref{C*eq} with $\eta=\Gr_3/2$, {\it i.e.},
\beq\label{Gt0=}
\Gt_0:= \frac {C\sqrt \Ge}{\frac{\Gr_3}{2} + \sqrt \Ge}.
\eeq
Define $\Phi := \Phi_3 \circ \Phi_2$, and $R_{\Gt_0} := (a_0, 0) \times (\Gt_0, 2\pi-\Gt_0)$. Then $\Phi^{-1}(R_{\Gt_0}) \cap \GO_1$ is a bounded subset of $\GO_1$. Define
\beq
\GO_{\Gt_0}:= \Phi(\Phi^{-1}(R_{\Gt_0}) \cap \GO_1).
\eeq
Then $\GO_{\Gt_0}$ is a connected subset of $R_{\Gt_0}$ and has two lateral boundaries denoted by $l_1$ and $l_2$.
Let
\beq
\check r_\Ge(r, \Gt) := (\tilde r_\Ge \circ \Phi^{-1})(r, \Gt) , \quad (r, \Gt) \in \GO_{\Gt_0},
\eeq where $\tilde r_\Ge$ is given in \eqref{tilderGe2}.
Then, $\check r_\Ge$ satisfies
\beq \left\{
\begin{array}{ll}
\ds\Delta \check r_{\Ge} = 0  &\quad\mbox{in } \GO_{\Gt_0} ,\\
\ds \check r_{\Ge} = 0  &\quad\mbox{on } l_1 \cup l_2 .
\end{array}
\right.
\eeq

We have the following lemma whose proof will be given at the end of this section.
\begin{lem}\label{lem:tr}
There is a constant $C$ such that for $(r, \Gt) \in \GO_{\Gt_0}$
\beq\label{tilderGe}
\left| \check r_\Ge (r,\theta)\right| \leq C  \exp  \left( - \frac{\pi}{|a_0|}  (\Gt -\Gt_0) \right) \quad\mbox{if } \Gt \leq \pi
\eeq
and
\beq\label{tilderGe3}
\left| \check r_\Ge (r,\theta)\right| \leq C  \exp  \left( - \frac{\pi}{|a_0|}  (2\pi-\Gt_0 -\Gt) \right) \quad\mbox{if } \Gt \geq \pi.
\eeq
\end{lem}

The desired inequality \eqnref{mainest3} follows from \eqnref{tilderGe} and \eqnref{tilderGe3}. To see this, we first observe that if $r+i\Gt= \Phi_3 \circ \Phi_2(w)$, then $e^{r+i\Gt}= \Phi_2(w)$, in other words, $\Gt=\mbox{arg}\Phi_2(w)$. Because of \eqnref{a=-Gg} and \eqnref{Gt0=}, we have
$\Gt_0/|a_0| \le C$ for some constant $C$ independent of $\Ge$ provided that $\Ge$ is sufficiently small. Observe that if $w=\xi+i\eta \in S$ and $\eta>0$, then $\Gt = \mbox{arg} \Phi_2(w) < \pi$. So it follows from Lemma \ref{C_*} and \eqnref{tilderGe} that
$$
|\tilde r_\Ge(w)| \le C \exp \left( - \frac{\pi}{|a_0|} \Gt \right) \le C_1 \exp \left(-  \frac {A}{\sqrt {\Ge} + |\eta|}\right).
$$ For $w=u+iv \in S$ with $\eta \leq 0$, $\Gt = \mbox{arg} \Phi_2(w) \geq \pi$. Lemmas \ref{C_*} and \ref{lem:tr} also yield
 $$
|\tilde r_\Ge(w)| \le C \exp \left( - \frac{\pi}{|a_0|} (2\pi - \Gt) \right) \le C_2 \exp \left(-  \frac {A_1}{\sqrt {\Ge} + |\eta|}\right).
$$
So we have \eqnref{mainest3} and the proof of Lemma \ref{d**} is completed.
\qed

\medskip

Let us now prove Lemma \ref{C_*} and Lemma \ref{lem:tr}.

\noindent{\sl Proof of Lemma \ref{C_*}}.  In this proof, we shall consider  the case when $w  = \xi + \eta i \in S$ with $\eta \ge 0$ only.  We first note that
$$
\mbox{Im}\,  \Phi_2 (w) = \frac{\eta(1-\Ga^2)}{(1-\Ga \xi)^2 + \Ga^2 \eta^2}.
$$
Using \eqnref{Ga=-1} one can see that
$$
1-\Ga^2 \ge C \sqrt{\Ge}
$$
if $\Ge$ is sufficiently small, since $|\xi|\leq \frac 1 2 \rho_3 \leq \frac 1 2 $.  We observe that for $w= \xi+i\eta \in S$,
$$ 1+\xi \le 1+c_3 -\sqrt{\rho_3^2 - \eta^2}  = 1+c_3-\rho_3 + \frac{\eta^2}{\rho_3 +\sqrt{\rho_3^2 - \eta^2} }  \le \Ge + \frac{\eta^2}{\rho_3}, $$ since $w \in B_4 \setminus B_3$ and $
\mbox{dist} \, (\p B_3, \p B_4 )  = - \frac{1}{1+\Ge} +1 \leq \Ge$ by \eqref{distB_3B_4}.

If $\eta \ge \sqrt{\Ge}$, then
$$\mbox{Im}\,  \Phi_2 (w) \ge C_1 \frac{\eta\sqrt{\Ge}}{\Ge + \eta^2} \ge C_2 \frac{\sqrt{\Ge}}{\eta}
$$ by \eqref {Ga=-1} and the property that $\frac{\eta}{\rho_3} \leq 1$.
Since $|\Phi_2(w)| \ge \frac 1 2 $ by \eqref{Gr_5=1}, we have
$$
\sin (\mbox{arg} \Phi_2 (w)) = \frac{\mbox{Im}\,  \Phi_2 (w)}{|\Phi_2(w)|} \ge 2 C_2 \frac{\sqrt{\Ge}}{\eta}.
$$
Thus we have
\beq\label{ArgPhi_2}
\mbox{arg} \Phi_2 (w) \ge C_3 \frac{\sqrt{\Ge}}{\eta}.
\eeq

If $0 \le \eta <\sqrt{\Ge}$, then there exists $w_0 = \xi_0 +  i\eta_0$ with $|\eta_0| =  \sqrt \Ge$ so that
$$
\mbox{arg} \Phi_2 (w) \geq \mbox{arg} \Phi_2 (w_0),
$$
so it follows from \eqnref{ArgPhi_2} that
$$
\mbox{arg} \Phi_2 (w) \ge C_3.
$$
This proves \eqnref{C*eq}.
\qed

\medskip

\noindent{\sl Proof of Lemma \ref{lem:tr}}. By definition, $\Omega_{\Gt_0}$ is a subset of $R_{\Gt_0}$, and $\p \Omega_{\Gt_0} \cap \p R_{\Gt_0}$ belongs to  $\Gt=\Gt_0$ or  $\ 2\pi-\Gt_0$. We define functions $\psi_\pm$ in $R_{\Gt_0}$ as the solutions to
\beq \left\{
\begin{array}{ll}
\ds\Delta \psi_{\pm}  = 0  \quad&\mbox{in } R_{\Gt_0} ,\\
\ds \psi_{\pm}(r, \Gt) = 0  \quad&\mbox{on } \p R_{\Gt_0} \setminus \p \Omega_{\Gt_0}  , \\
\ds \psi_{\pm}(r, \Gt) =\max \left\{\pm\check r_{\Ge}(r, \Gt) ,0\right\} \quad&\mbox{on } \p R_{\Gt_0} \cap \p \Omega_{\Gt_0} .
\end{array}
\right.
\eeq It was shown in \cite{ACKLY} that $\norm{r_{\Ge}-h}_{L^{\infty} (\GO)}$ is bounded independently of $\Ge$. So, there is a constant $M$ independent of $\Ge>0$ such that
\beq\label{breveM}
|\psi_{\pm}(r, \Gt)| \leq M \quad\mbox{for all } (r, \Gt) \in R_{\Gt_0},
\eeq
and it can be shown in the same way as \eqref{V-UV+} in the previous section that
\beq\label{R_rR+}
-\psi_- \leq   \check r_{\Ge} \leq \psi_+  \quad\mbox{in } \GO_{\Gt_0}.
\eeq
So to prove \eqnref{tilderGe3} it suffices to show
\beq\label{psipm}
\left| \psi_{\pm} (r,\theta)\right| \leq C  \exp  \left( - \frac{\pi}{|a_0|}  (2\pi-\Gt_0 -\Gt) \right)  ~\mbox{for } \Gt \in [\pi,2\pi-\theta_0].
\eeq

We prove \eqnref{psipm} only for $\psi_+$ since the proof for $\psi_-$ is identical. The solution $\psi_+$ can be found explicitly:
$$
\psi_+ = \psi_+^e + \psi_+^o
$$
where $\psi_+^e$ and $\psi_+^o$ are the even and odd parts about $ \theta = \pi$ given by
\begin{align*}
\psi_+^e (r,\Gt) &= \sum_{n=1} ^{\infty} \Ga_n \sin \left( \frac {n\pi} {|a_0|} r  \right)
\cosh \left( \frac { n \pi } {|a_0|}\left (\Gt-  \pi  \right) \right)   \\
\psi_+^o (r,\Gt) &= \sum_{n=1} ^{\infty} \Gb_n \sin \left( \frac {n\pi} {|a_0|} r  \right) \sinh  \left( \frac { n \pi } {|a_0|}\left (\Gt- \pi  \right) \right)
\end{align*}
for some constants $\Ga_n$ and $\Gb_n$.

Suppose that $\Gt \geq \pi$. Then we have
\begin{align*}
\left| \psi_+^e (r,\Gt) \right| & \leq  \sum_{n=1}^{\infty}   | \Ga_n |  \exp \left( \frac{n \pi}{|a_0|} (\Gt-  \pi) \right) \\
& \leq 2 \sum_{n=1} ^{\infty}   | \Ga_n | \cosh  \left( \frac{n \pi}{|a_0|} (\pi -\Gt_0) \right)
\exp \left( -\frac{n \pi}{|a_0|} (2\pi- \Gt_0 -\Gt) \right)
\end{align*}
Note that
$$
\left( \sum_{n=1}^{\infty}  | \Ga_n |^2 \cosh^2 \left( \frac{n \pi}{|a_0|} (\pi -\Gt_0) \right) \right)^{1/2} =\left(\frac 2 {|a_0|}\right)^{1/2} \norm{   \psi_{+}^e (\cdot,2\pi - \Gt_0)  }_{L^2 ([a_0,0])} \le \sqrt 2 M,
$$ since $\norm{   \sin \left( \frac {n\pi} {|a_0|} r  \right)  }_{L^2 ([a_0,0])} = \left(\frac 1 2 |a_0|\right)^{1/2} $.
So it follows from the Cauchy-Schwarz inequality that
\begin{align*}
\left| \psi_+^e (r,\Gt) \right| \le 2\sqrt 2M  \left(\sum_{n=1}^{\infty} \exp \left( -\frac{2n \pi}{|a_0|} (2\pi- \Gt_0 -\Gt) \right) \right)^{1/2},
\end{align*}
and hence
\beq\label{psi+e}
\left| \psi_+^e (r,\Gt) \right| \le C \exp \left( -\frac{\pi}{|a_0|} (2\pi- \Gt_0 -\Gt) \right)
\eeq
for some constant $C$. Since $\psi_+^e(r, \Gt) = \psi_+^e(r, 2\pi- \Gt)$, $$
\left| \psi_+^e (r,\Gt) \right| \le C \exp \left( -\frac{\pi}{|a_0|} (\Gt- \Gt_0 ) \right)
$$ when $\Gt < \pi$ as well.

Since  $\sinh B \leq \sinh A \, e^{B-A}$ if $0 < B < A$, we obtain, for $\Gt \geq \pi  $,
\begin{align*}
\left| \psi_{+}^o (r,\Gt)\right| & \leq  \sum_{n=1}^{\infty} | \Gb_n |  \sinh \left( \frac{n \pi}{|a_0|} \left (\Gt- \pi  \right) \right) \\
& \leq  \sum_{n=1}^{\infty} | \Gb_n |  \sinh \left( \frac {n \pi}{|a_0|} (\pi-\Gt_0) \right)
\exp  \left( -\frac {n \pi}{|a_0|} (2\pi -\Gt_0- \Gt)  \right),
\end{align*}
and hence
\beq\label{psi+o}
\left| \psi_+^o (r,\Gt) \right| \le C \exp \left( -\frac{\pi}{|a_0|} (2\pi- \Gt_0 -\Gt) \right).
\eeq
Because of symmetry of $\psi_+^o$, we have
$$
\left| \psi_+^o (r,\Gt) \right| \le C \exp \left( -\frac{\pi}{|a_0|} (\Gt- \Gt_0 ) \right),
$$
when $\Gt < \pi$ as well.  So we have \eqnref{psipm} and the proof is complete. \qed

%%%%%%%%%%%%%%%%%%%%%%%%%%%%%%%%%%%%%%%%%%%%%%%%%%%%%%%%%%
\section{Proofs of Theorem \ref{1st_main}}\label{sec:proof}
%%%%%%%%%%%%%%%%%%%%%%%%%%%%%%%%%%%%%%%%%%%%%%%%%%%%%%%%%%

In this section we prove Theorem \ref{1st_main} and  \eqnref{limit}.

\noindent {\sl Proof of Theorem \ref{1st_main}}. \ One can see from \eqnref{intzero}, \eqnref{q_equ} and \eqnref{decom1} that
\beq \label{GaGe_int}
\Ga_{\Ge} = \int_{\p D_1^0}  \p_{\nu} r_\Ge~ds .
\eeq
So, it is enough to prove
\beq\label{ineqaulity1}
\left|\int_{\p D_1^0} \p_{\nu} ( r_{\Ge} -  u_{0}) ds \right| \leq C {\Ge} |\log \Ge|^{2m-1}
\eeq
for some constant $C$ independent of $\Ge$.

Let $V:= \Rbb^2 \setminus \overline{(D_1^0 \cup D_2^0 \cup D_2^\Ge)}$, and let $\GG_1:= \p D_2^0 \setminus D_2^\Ge$ and $\GG_2: =\p D_2^\Ge \setminus D_2^0$. Then, $\p D_1$, $\GG_1$, and $\GG_2$ constitute the boundary of $V$. Let
\beq
\Gvf_\Ge (\Bx):= r_\Ge(\Bx) - u_0(\Bx) - (r_\Ge(0,0)- u_0(0,0)).
\eeq
Then, $\Gvf_\Ge$ is a bounded harmonic function in $V$ and $\Gvf_\Ge \equiv 0$ on $\p D_1$. We claim that
\beq\label{ieq_crucial_2}
| \Gvf_\Ge (\Bx) | \leq C \Ge, \quad \Bx \in V.
\eeq
In fact, if $\Bx \in \GG_1$, then $u_0(\Bx) - u_0(0,0)=0$ and $r_\Ge(\Bx+\Ge) - r_\Ge(0,0)=0$. Therefore, since $\nabla r_\Ge$ is bounded on any bounded subset of $\mathbb{R}^2 \setminus \overline{D_1 \cup D_2}$ (refer to \cite{ACKLY}), we have
\beq\label{bound_b}
|\Gvf_\Ge (\Bx)| = |r_\Ge(\Bx) - r_\Ge(0,0)|
= |r_\Ge(\Bx) - r_\Ge(\Bx+\Ge)| \le C \Ge.
\eeq
Likewise we have for $\Bx \in \GG_2$
\beq\label{bound_b2}
|\Gvf_\Ge (\Bx)| \le C \Ge.
\eeq
Since $\Gvf_\Ge \equiv 0$ on $\p D_1$ and $\Gvf_\Ge (\Bx)$ is bounded, we obtain \eqnref{ieq_crucial_2} by the maximum principle.

Choose $M$ so large that
$$
D_1^0 \subset \left(-\frac M 2, 0 \right)\times\left(-\frac M 2, \frac M 2\right),
$$
and let $\Go=(-M, 0)\times(-M,M)$. Since $\Gvf_\Ge$ is harmonic in $V$, we have
\begin{align*}
\int_{\p D_1^0} \p_{\nu} ( r_{\Ge} -  u_{0}) ds  = \int_{\p D_1^0} \p_{\nu} \Gvf_\Ge ds = \int_{\p \Go} \p_{\nu} \Gvf_\Ge ds.
\end{align*}

Divide $\p\Go$ into three pieces: $\p\Go=\Gg_1 \cup \Gg_2 \cup \Gg_3$ where
$$
\Gg_1 := \left\{ (0,y)~|~ |y| \le \frac {A_0} {|\log {\Ge}|}  \right\}, \quad \Gg_2 :=\left \{ (0,y)~|~ \frac {A_0} {|\log {\Ge}|}  < |y| \le M \right \},\quad \Gg_3:= \p\Go \setminus (\Gg_1 \cup \Gg_2),
$$
and write
$$
\int_{\p \Go} \p_{\nu} \Gvf_\Ge ds = \int_{\Gg_1} + \int_{\Gg_2} + \int_{\Gg_3}  \p_{\nu} \Gvf_\Ge ds := I + II + III.
$$ Here, the constant $A_0$ is given by Theorems \ref{asymdecay} and \ref{asymdecay2} so that \begin{equation} |\nabla u_0(x,y)| + |\nabla r_0(x,y) | \leq C \exp \left( -\frac {A_0} {|y| + \sqrt {\Ge}}\right) \label{gamma_1} \end{equation}
for $|y|< \Gd$ and $x \in (x_1(y),x_2(y))$.

If $-\frac {A_0}{|\log {\Ge}|} \le y \le\frac {A_0} {|\log {\Ge}|}  $, then \eqref{gamma_1} implies that
$$
|\nabla \Gvf_\Ge (0,y)|  \leq C \exp \left( -\frac {A_0} {|y| + \sqrt {\Ge}}\right). $$
Thus we have
\beq\label{I_esti}
|I| \le C {\Ge}.
\eeq

If $\frac {A_0}{|\log {\Ge}|} < |y| \le M$, there is $r > Cy^{2m}$ for some $C$ such that $B_r(0,y) \subset V$. It then follows from a gradient estimate for harmonic functions and \eqnref{ieq_crucial_2} that
$$
\left|\nabla \Gvf_\Ge (0,y)\right| \leq C \frac {\Ge}{y^{2m}} ,
$$
and
\beq\label{II_esti}
|II| \le C \Ge \int_{\frac {A_0}{|\log {\Ge}|} < |y| \le M} \frac{1}{y^{2m}} dy \leq  C {\Ge} |\log \Ge|^{2m-1}.
\eeq

There is a constant $r>0$ such that $B_r(\Bx) \subset V$ for all $\Bx \in \Gg_3$. So, we have
from \eqref{ieq_crucial_2} that for any $\Bx \in \Gg_3$,
$$
\left|\nabla \Gvf_\Ge (\Bx)\right| \leq C \frac {\Ge}{r} \leq C \Ge,
$$
and
\beq\label{III_esti}
|III| \leq  C \Ge.
\eeq
Now, \eqref{ineqaulity1} follows from \eqref{I_esti}, \eqref{II_esti}, and \eqref{III_esti}, and the proof is complete.
\qed

\medskip
The formula \eqnref{limit} is an immediate consequence of \eqnref{mainsing}. In fact, if  $r_1$ and $r_2$ are radii of circles osculating to $\p D_1$ and $\p D_2$ at $(0,0)$ and $(\Ge,0)$, respectively, then it is proved in \cite{LY} that $\Bp_1$ and $\Bp_2$ which are fixed points of mixed reflections are given by
\beq
\Bp_i = \left( (-1)^i \sqrt 2 \sqrt {\frac {r_1 r_2}{r_1 + r_2}}\sqrt {\Ge} + O (\Ge), 0\right) \quad\mbox{as } \Ge \rightarrow 0. \label{p_i_esti}
\eeq
So we obtain \eqnref{limit} from \eqnref{mainsing}.

%%%%%%%%%%%%%%%%%%%%%%%%%%%%%%%%%%%%%%%%%%%%%%%%%%%%
\section{Approximations  of $\Ga_0$}\label{sec:approx}
%%%%%%%%%%%%%%%%%%%%%%%%%%%%%%%%%%%%%%%%%%%%%%%%%%%%

The region outside $D_1^0 \cup D_2^0$ has cusps at $(0,0)$, and it may cause some problem in computing $\Ga_0$. To avoid this trouble, we show that by replacing the cusp with a neck a good approximation of $\Ga_0$ can be obtained.

For $\rho>0$ let
\beq
D_{(\rho)} = \left( D_1^0 \cup D_2^0 \right)  \cup  \left([-\rho,\rho] \times  [-\rho,\rho]\right)
\eeq
which is of dumbbell shape, and let $u_{(\rho)}$ be the solution to
\beq
\left\{
\begin{array}{ll}
\ds \Delta u_{(\rho)}  = 0 \quad\mbox{in } \Rbb^2 \setminus \overline{D_{(\rho)}},\\
\ds u_{(\rho)} = \Gl_{(\rho)} \ (\mbox{constant}) \quad\mbox{on } \p D_{(\rho)}, \\
\ds u_{(\rho)}(\Bx) - h(\Bx) =O(|\Bx|^{-1}) \quad \mbox{as }
|\Bx|\rightarrow \infty,
\end{array}
\right.
\eeq
where the constant $\Gl_{(\rho)}$ is determined by the additional condition
\beq
 \int_{\p  D_{(\rho)} } \frac {\p u_{(\rho)}}{\p \nu} \Big|_+  ds  =0.
\eeq

We have the following theorem.
\begin{thm}\label{2nd_main}
Let $\Gd$ be the number appearing in Theorem \ref{thmtouching}. For $\rho \in (0, \Gd/2)$,  let
\beq
\Ga_{(\rho)} = \int_{\p D_1^0 \setminus  [-2\rho,2\rho] \times  [-2\rho,2\rho]  } \p_{\nu} u_{(\rho)}~ ds.
\eeq
Then, there are constants $C$ and $A$ such that
\beq
\left| \Ga_0 - \Ga_{(\rho)}\right| \leq C \exp \left( - \frac A \rho \right).
\eeq
\end{thm}

\pf
Choose a point $z_0$  on the common boundary of $D_{(\Gr)}$ and $D_1^0 \cup D_2^0$ and let
$$
\Gvf(z):= u_{(\rho)}(z) -u_0(z)  - (u_{(\rho)}  (z_0)-u_0  (z_0)).
$$
Since $u_{(\rho)}(z) -  u_{(\rho)} (z_0)=0$ for all $z \in \p D_{(\Gr)}$ and $u_{0}(z)  - u_{0} (z_0) =0$ on $\p D_1^0 \cup \p D_2^0$, we have
\beq\label{estone}
\Gvf(z)=0, \quad z \in \p D_{(\rho)} \setminus  \left( [-\rho,\rho] \times  [-\rho,\rho] \right).
\eeq
On the other hand, if $x_1(\rho)\le x \le x_2 (\rho)$, then we have from \eqnref{u0decay}
$$
| u_{0} (x + i\rho) -  u_{0} (z_0) | \le C e^{- \frac{A}{\rho}},
$$
and hence
\beq\label{esttwo}
| \Gvf (x+i\rho) |= |u_0 (x,\rho)-u_{0} (z_0)| \le C e^{- \frac{A}{\rho}}.
\eeq
Similarly one can see that if $x_1(-\rho)\le x \le x_2 (-\rho)$, then
\beq\label{estthree}
| \Gvf (x-i\rho) | \le
C e^{- \frac{A}{\rho}}.
\eeq
It follows from \eqnref{estone}, \eqnref{esttwo}, and \eqnref{estthree} that
\beq\label{estfour}
| \Gvf (z) | \le C e^{- \frac{A}{\rho}}
\eeq
for all $z \in \p D_{(\rho)}$, and hence for all $z \in \mathbb{R}^2 \setminus \overline{D_{(\rho)}}$ by the maximum principle. Note that we may apply the maximum principle since $u_{(\Gr)}(z)-u_0(z) \to 0$ as $|z| \to \infty$.

We now estimate $\nabla (u_{(\rho)}  (z)-u_0  (z))= \nabla \Gvf(z)$ on $\p D_1^0 \setminus  [-2\rho,2\rho] \times  [-2\rho,2\rho]$. Because of \eqnref{estone}, one can apply the argument used right after of Theorem \ref{thmtouching} to see that $\Gvf(z)$ can be extended across $\p D_1^0 \setminus  [-2\rho,2\rho] \times  [-2\rho,2\rho]$ so that the extended function is harmonic in $\overline{B_r(z)}$ for all $z \in \p D_1^0 \setminus  [-2\rho,2\rho] \times  [-2\rho,2\rho]$ where $r= s \Gr^{2m}$ for some $s>0$ (independent of $\Gr$). Then by the gradient estimate for harmonic functions we have
\beq\label{estfive}
|\nabla (u_{(\rho)} -u_0) (z)| \le C_2  e^{- \frac{A_2}{\rho}}, \quad z \in \p D_1^0 \setminus  [-2\rho,2\rho] \times  [-2\rho,2\rho].
\eeq

It then follows from \eqnref{asymdecay} and \eqnref{estfive} that
\begin{align*}
|  \Ga_0 - \Ga_{(\rho)} | &\leq \int_{(\p D_1^0) \setminus  [-2\rho,2\rho] \times  [-2\rho,2\rho]  } \left| \p_{\nu} ( u_{(\rho)} - u_0) \right|\, ds + \int_{(\p D_1^0) \cap (  [-2\rho,2\rho] \times  [-2\rho,2\rho])  } \left| \p_{\nu} u_{0}\right|\, ds \\
& \leq C {e^{- \frac {A_3}{\rho} }}.
\end{align*}
This completes the proof.
\qed

%%%%%%%%%%%%%%%%%%%%%%%%%%%%%%%%%%%%%%%%%%%%%%%%%%%%%%%%


\begin{thebibliography}{10}

\bibitem{ACKLY} H. Ammari, G. Ciraolo, H. Kang, H. Lee and K. Yun, Spectral analysis of the Neumann-Poincar\'e operator and characterization of the stress concentration in anti-plane elasticity, Arch. Ration. Mech. An. 208 (2013), 275--304.

\bibitem{AKLLL} { H. Ammari, H. Kang, H. Lee, J. Lee and M. Lim}, { Optimal bounds on the gradient of solutions to conductivity problems}, J. Math. Pure. Appl. 88 (2007), 307--324.

\bibitem{AKLLZ} H. Ammari, H. Kang, H. Lee, M. Lim and H. Zribi, Decomposition theorems and fine estimates for electrical fields in the presence of closely located circular inclusions, J. Differ. Equations 247 (2009), 2897-2912.

\bibitem{AKLim} {H. Ammari, H. Kang and M. Lim}, {Gradient estimates for solutions to the conductivity problem}, Math. Ann. 332(2) (2005), 277--286.

\bibitem{bab} I. Babu\u{s}ka, B. Andersson, P. Smith and K. Levin, Damage analysis of fiber composites. I. Statistical
analysis on fiber scale, Comput. Methods Appl. Mech. Engrg. 172 (1999), 27--77.

\bibitem{BLY} E. Bao, Y.Y. Li, B. Yin, Gradient estimates for the perfect conductivity problem, Arch. Ration. Mech. An. 193 (2009), 195-226.

\bibitem{BLY2} {E. S. Bao, Y.Y. Li and B. Yin}, {Gradient estimates for the perfect and insulated conductivity problems with multiple inclusions}, Commun. Part. Diff. Eq. 35 (2010), 1982--2006.

\bibitem{BKN} L. Berlyand, A. Kolpakov and A. Novikov, {\sl Introduction to the Network Approximation Method for Materials Modeling}, Encyclopedia of Mathematics and its Applications, Cambridge, Cambridge Univ. Press, 2012.

\bibitem{bv} {E. Bonnetier and  M. Vogelius}, {An elliptic regularity result for a composite medium with ``touching" fibers of circular cross-section}, SIAM J. Math. Anal. 31 (2000), 651--677.

\bibitem{BC} {B. Budiansky and G. F. Carrier}, { High shear stresses in stiff fiber composites}, J. Appl. Mech. 51 (1984),  733--735.

\bibitem{CG} H. Cheng and L. Greengard,
\newblock A method of images for the evaluation of electrostatic
fields in systems of closely spaced conducting cylinders,
\newblock SIAM J. Appl. Math., 58 (1998), 122--141.

\bibitem{KLY} H. Kang, M. Lim and K. Yun, Asymptotics and computation of the solution to the conductivity
equation in the presence of adjacent inclusions with extreme
conductivities, J. Math. Pure. Appl. 99 (2013), 234--249.

\bibitem{KLY13} H. Kang, M. Lim and K. Yun, Characterization of the electric field concentration between two adjacent spherical perfect conductors, SIAM J. Appl. Math., to appear.

\bibitem{keller} J.B. Keller, Stresses in narrow regions, Trans. ASME J. Appl. Mech., 60 (1993), 1054--1056.

\bibitem{lekner10} J. Lekner, Analytical expression for the electric field enhancement between two closely-spaced conducting spheres, J. Electrostatics 68 (2010), 299-304.

\bibitem{lekner11} J. Lekner, Near approach of two conducting spheres: enhancement of external electric field, J. Electrostatics 69 (2011), 559-563.

\bibitem{lekner} J. Lekner, Electrostatics of two charged conducting spheres, Proc. R. Soc. A, 468 (2012), pp. 2829-2848.

\bibitem{LLBY} H. Li, Y.Y. Li, E.S. Bao and B. Yin, Derivative estimates of solutions of elliptic systems in narrow regions, Quart. Appl. Math., to apear.

\bibitem{ln} {Y.Y. Li and L. Nirenberg}, { Estimates for elliptic system from composite material},
Comm. Pure Appl. Math., LVI (2003), 892--925.

\bibitem{lv} { Y.Y. Li and M. Vogelius}, {Gradient estimates for solution to divergence form elliptic equation with discontinuous coefficients}, Arch. Rat. Mech. Anal. 153 (2000), 91--151.

\bibitem{LY} {M. Lim and K. Yun}, {Blow-up of electric fields between closely spaced spherical perfect conductors},
Commun. Part. Diff. Eq. 34 (2009), 1287--1315.

\bibitem{LY2} {M. Lim and K. Yun}, {Strong influence of a small fiber on shear stress in fiber-reinforced composites},
J. Differ. Equations 250 (2011), 2402--2439.

\bibitem{MPM} R.C. McPhedran, L. Poladian and G.W. Milton, Asymptotic studies of closely spaced, highly conducting cylinders, Proc. R. Soc. Lond. A 415 (1988), 185--196.

\bibitem{MS} V. G. Maz'ya and A.A. Solov'ev, On an integral equation for the Dirichlet problem in a plane domain with cusps on the boundary, {Math. USSR Sb.} 68 (1991), 61--83.

\bibitem{milton_book} G. W. Milton, {\sl  The Theory of Composites},
\newblock Cambridge Monographs on Applied and Computational
Mathematics, Cambridge, Cambridge Univ. Press, 2001.

\bibitem{Y} {K. Yun},  {Estimates for electric fields blown up between closely adjacent conductors with arbitrary shape}, {SIAM J. Appl. Math.} 67 (2007), 714--730.

\bibitem{Y2} {K. Yun}, {Optimal bound on high stresses occurring between stiff fibers with arbitrary shaped cross sections}, J. Math. Anal. Appl. 350 (2009), 306-312.

\end{thebibliography}
\end{document}